\documentclass[12pt]{article}
\usepackage{group21}

\title{Reversing Symmetry Groups of Cat Maps}

\author{Michael Baake}
\address{Institut f\"ur Theoretische Physik, Universit\"at T\"ubingen,\\
         Auf der Morgenstelle 14, D-72076 T\"ubingen, Germany}
\newcommand{\zr}[1]{\mbox{\hspace*{#1em}}}

\newcommand{\ZZ}{\mbox{\sf Z\zr{-0.45}Z}}

\newcommand{\Id}{\mbox{\rm 1\zr{-0.64}{\small 1}}}
\newcommand{\TT}{\mbox{\sf T}}

\newtheorem{lemma}{Lemma}
\newtheorem{theorem}{Theorem}

\begin{document}
\maketitle
\begin{abstract}
Toral automorphisms are widely used (discrete) dynamical systems,
the perhaps most prominent example (in 2D) being Arnold's cat map.
Given such an automorphism $M$, its {\em symmetries} (i.e.\ all
automorphisms that commute with $M$) and {\em reversing symmetries}
(i.e.\ all automorphisms that conjugate $M$ into its inverse)
can be determined by means of number theoretic tools. 

Here, the case of $Gl(2,\ZZ)$ is presented and the possible 
(reversing) symmetry groups are completely classified. 
Extensions to affine mappings and to 
$k$-(reversing) symmetries (i.e.\ (reversing) symmetries of $M^k$),
and applications to the projective group $PGl(2,\ZZ)$ and to 
trace maps are briefly discussed.
\end{abstract}

\section*{Introduction}

Symmetries of a dynamical system provide useful tools to understand its
phase space structure and to simplify its description. In what follows,
we investigate, in a group theoretic setting, symmetries of toral automorphisms,
i.e.\ of matrices $M \in Gl(2,\ZZ)$. This is similar in spirit to several
recent publications \cite{BR1,BR2,Lamb} where the (discrete) dynamical 
system, $F$, is an element of some automorphism group, ${\cal G}$ 
($Gl(2,\ZZ)$ in our case).
Then, the group of {\em symmetries} of $F$ is the centralizer of $F$ in
${\cal G}$, defined as
\begin{equation}
    {\cal S} (F) \; := \; \{ G \in {\cal G} \mid F\circ G = G\circ F \} \, .
\end{equation}

Beyond a symmetry, one could also have a so-called {\em reversing symmetry}
which maps $F$ to its inverse, $G \circ F \circ G^{-1} \; = \; F^{-1}$.
All such reversing symmetries, together with the symmetries, form another
group, the {\em reversing symmetry group} of $F$,
\begin{equation}
    {\cal R} (F) \; := \; \{ G \in {\cal G} \mid 
                             G \circ F \circ G^{-1} = F^{\pm 1} \} \, .
\end{equation}
If $F$ is an involution (i.e.\ $F^2=\Id\,$) or if $F$ admits no reversing
symmetry, one has ${\cal R} (F) = {\cal S} (F)$. Otherwise, ${\cal R} (F)$
contains ${\cal S} (F)$ as a normal subgroup of index 2.

One cannot expect to get a full characterization of the (reversing) 
symmetry group of an arbitrary automorphism $F$. In special situations,
however, a classification is possible. The latter is only necessary
mod conjugacy, as ${\cal R} (H \circ F \circ H^{-1}) = 
H \circ {\cal R} (F) \circ H^{-1}$. In our present context, for the cat maps,
we are in the lucky situation of a matrix group, ${\cal G} = Gl(2,\ZZ)$.
Since ${\cal G}$ is not continuous, it is not surprising that methods
from number theory are extremely handy to solve the puzzle.

\section*{Results: Symmetries}
Let us first investigate the structure of the symmetry group of a unimodular
$2\!\times\!2$-matrix. There are three types of matrices, called {\em elliptic}
(i.e.\ of finite order), {\em parabolic} (i.e.\ eigenvalues $\pm 1$, but not
of finite order), and {\em hyperbolic} (all remaining cases). The latter ones
are the {\em cat maps}. {}For the sake of completeness, the results will be
given for {\em all} types of matrices, not just for the cat maps.

Before we actually state the result, let us briefly sketch how to get there.
The case of elliptic and parabolic elements is most easily tackled by
finding the corresponding conjugacy classes and by determining the
corresponding symmetry groups for suitable representatives explicitly.
Then, the hyperbolic matrices have the property that their two eigenvalues
are irrational and different from one another. One can then employ Dirichlet's
unit theorem from algebraic number theory, specialized to the case of quadratic
number fields, to determine the corresponding unit groups and their relation
to the symmetry groups, for details we refer to \cite{BR2}.
The result is:
\begin{theorem} \label{thm1}
The structure of the centralizer of an element $M \in Gl(2,\ZZ)$,
i.e.\ the structure of the symmetry group ${\cal S}(M)$,
has precisely one of the following forms: \\
1) ${\cal S}(M) = Gl(2,\ZZ)$ if and only if $M=\pm \, \Id\,$, \\
2) ${\cal S}(M) = \{\pm \, \Id\,, \pm M \} \simeq C_2 \times C_2$
    if and only if ${\rm tr}(M)=0$ and $\det(M)=-1$ \\
    \hspace*{1em} (i.e.\ $M^2 = \Id \neq \pm M$), \\
3) ${\cal S}(M) = \{\pm \, \Id\,, \pm M \} \simeq C_4$
    if and only if ${\rm tr}(M)=0$ and $\det(M)=1$ \\
     \hspace*{1em} (i.e.\ $M^4=-M^2 = \Id\,$), \\
4) ${\cal S}(M) = \{\pm \, \Id\,, \pm M, \pm M^2 \} \simeq C_6$
    if and only if ${\rm tr}(M)=\pm 1$ and $\det(M)=1$ \\
     \hspace*{1em} (i.e.\ $M^6 = \Id \neq M^2$), or \\
5) ${\cal S}(M) \simeq C_2 \times C_{\infty}$
    if and only if $M$ is not of finite order.
\end{theorem} 

Case (5) needs some comment. If $M$ is not of finite order, it generates
an infinite cyclic group of symmetries (because $M$ commutes with any power
of itself), and it trivially commutes with $-\Id\,$. This also gives a group
of type $C_2 \times C_{\infty}$, but it need {\em not} be ${\cal S}(M)$
because there could be a {\em root} of $M$ still inside $Gl(2,\ZZ)$. The
non-trivial part of statement (5) is that, even if this happens, ${\cal S}(M)$
is still of type $C_2 \times C_{\infty}$ -- a direct consequence of Dirichlet's
unit theorem which also results in an explicit algorithm how to determine
the generator of $C_{\infty}$, see \cite{BR2} for details.

\section*{Results: Reversing Symmetries}
At this point, we move on to the {\em reversing} symmetries. Here, it turns
out that all elliptic and all parabolic elements are reversible, i.e.\ possess
a reversing symmetry (this is trivial for involutions, being their own reversing
symmetry). This can again most easily be checked explicitly for the conjugacy
class representatives. With the cat maps, however, things are more complicated.
It is {\em not} true that all cat maps are reversible, and it is, in general,
not easy to decide on reversibility versus irreversibility. Before we come
back to this point, let us, for the moment, just distinguish these two cases
and classify the possible reversing symmetry groups. If a matrix $M$ is
reversible, it can happen that there is a reversing symmetry of order 2,
or that the smallest order of a reversing symmetry is 4 -- resulting in
different groups, ${\cal R}(M)$. We obtain:
\begin{theorem} \label{thm1a}
The structure of the reversing symmetry group ${\cal R}(M) \subset
Gl(2,\ZZ)$ of a matrix $M \in Gl(2,\ZZ)$ has precisely one of the 
following forms: \\
1) ${\cal R}(M) = Gl(2,\ZZ)$ if and only if $M=\pm \, \Id\,$, \\
2) ${\cal R}(M) \simeq D_2$\,
    if and only if ${\rm tr}(M)=0$ and $\det(M)=-1$ 
    (i.e.\ $M^2 = \Id \neq \pm M$), \\
3) ${\cal R}(M) \simeq D_4$
    if and only if ${\rm tr}(M)=0$ and $\det(M)=1$ 
    (i.e.\ $M^4=-M^2 = \Id\,$), \\
4) ${\cal R}(M) \simeq D_6$
    if and only if ${\rm tr}(M)=\pm 1$ and $\det(M)=1$ 
    (i.e.\ $M^6 = \Id \neq M^2$), \\
5) ${\cal R}(M) \simeq D_{\infty} \times C_2$
    if and only if $M$ is of infinite order and possesses a \\
    \hspace*{1em} reversing symmetry of order 2, \\
6) ${\cal R}(M) \simeq C_{\infty} \times_s C_4$
    if and only if $M$ is of infinite order and possesses a \\
    \hspace*{1em} reversing symmetry of order 4, but none of order 2, or \\
7) ${\cal R}(M) \simeq C_{\infty} \times C_2$
    if and only if $M$ is of infinite order but irreversible.
\end{theorem}
Here, $D_n \simeq C_n \times_s C_2^{}$ is the dihedral group.
Let us add a comment on the structure of the reversing symmetry group
in case (6) of the last Theorem. It may look a bit astonishing that $\cal R$
can still be written as a semi-direct product, but the reason for it is that
the fourth order reversing symmetry $G$ fulfils $G^2=-\Id\,$, so by absorbing
the $C_2$-part of the symmetry group $\cal S$ we can find a subgroup
of $\cal R$ isomorphic to $C_4$ that conjugates the $C_{\infty}$-part
into itself but has only the unit matrix in common with it.

The next step of refinement would consist in a complete characterization
of those cat maps that are reversible. There are two possibilities to
proceed. One uses the fact that 
$Gl(2,\ZZ) \simeq (PSl(2,\ZZ) \times_s C_2) \times \{\pm\,\Id\,\}$.
Then, reversibility of $M$ in
$Gl(2,\ZZ)$ (i.e.\ $M$ conjugate to $M^{-1}$) can be reduced to various
relations between elements of $PSl(2,\ZZ)$. The latter can be decided upon
in finitely many steps because $PSl(2,\ZZ) \simeq C_2 * C_3$ (free product of
two cyclic groups) and the word problem is completely solvable there, see
\cite{BR2} for details. Another possibility uses a reformulation of the
equation for reversibility in terms of an integer quadratic form and results
in the question whether a certain integer can be represented by that form 
\cite{BR2}. This can then be decided by a finite algorithm which is the best 
type of result one may expect here.

\section*{Extensions}
So far, we have analyzed (reversing) symmetries of $M \in Gl(2,\ZZ)$.
A power of $M$ could, in principle, have additional (reversing) symmetries,
called {\em (reversing) $k$-symmetries} (if valid for $M^k$, but no smaller 
power of $M$). The answer, however, is essentially negative. 
In fact, we have (for $k>1$):
\begin{theorem}
 Elliptic elements of $Gl(2,\ZZ)$ cannot have (reversing) $k$-symmetries, 
 unless $M^k=\pm\,\Id\,$, where ${\cal R}(M^k) = {\cal S}(M^k) = Gl(2,\ZZ)$.
 Parabolic elements cannot have any (reversing) $k$-symmetries.
 {}Finally, hyperbolic elements cannot have $k$-sym\-me\-tries, and at most
 reversing 2-symmetries (if $\det(M)=-1$, but $M^2$ reversible).
\end{theorem}

The finite order and parabolic cases again rest upon the representatives of
the conjugacy classes, while the rest follows from special properties of
$2\!\times\!2$-matrices and the Cayley-Hamilton theorem for them.

Another rather obvious extension concerns the class of transformations that one
admits as symmetries. In particular, one could search for (reversing)
symmetries of a toral automorphism within the larger group of {\em affine}
transformations on the torus $\TT = [0,1)^2$, giving the semi-direct product 
${\cal G}_a = \TT\times_s Gl(2,\ZZ)$. This, indeed, gives rise to a number
of interesting possibilities because of the following result (where we write
the affine transformations as $(t,G)$ with obvious meaning).
\begin{lemma}
 The affine transformation $(t,G)$ is a (reversing) symmetry of the toral
 automorphism $(0,M)$ if and only if $G$ is a (reversing) symmetry of $M$
 in $Gl(2,\ZZ)$ and $Mt=t$ (mod $\TT$).
\end{lemma}

Let us remark that this also gives rise to (reversing) $k$-symmetries because
$M^k t=t$ has an increasing number of solutions on $\TT$. They can be counted,
provided no eigenvalue of $M^k$ is 1, as $a_k = |\det(M^k-\Id\,)|$. Grouping
them into orbits under the action of $M$, one also relates this symmetry problem
to the structure of dynamical or Artin-Mazur $\zeta$-functions which might be
an interesting side-remark.

\section*{Trace Maps}
The above results also admit the full treatment of the corresponding 
problem for matrices in the projective linear group 
$PGl(2,\ZZ) = Gl(2,\ZZ)/\{\pm\,\Id\,\}$.
One can read off the even simpler classification of possible (reversing)
symmetry groups from Theorems 1 and 2, rederiving a result \cite{BR1}
previously obtained directly (though with similar techniques). 
What is more, Theorem 3
has a counterpart where now even the possibility of reversing 2-symmetries
vanishes, because all equations are now ``mod $\pm\,\Id\,$'' which makes
many orientation reversing matrices actually reversible, compare \cite{BR2}.

The interest in this case originates from the isomorphism between $PGl(2,\ZZ)$ 
and the group of so-called Nielsen trace maps \cite{BGJ}. They are the 
(invertible) polynomial mappings of 3-space into itself that preserve the 
{}Fricke-$\!$Vogt invariant
\begin{equation}
   I(x,y,z) \; = \; x^2 + y^2 + z^2 - 2 x y z - 1
\end{equation}
and fix the point $(1,1,1)$, the Fibonacci trace map being its
best-studied example.

The surface $\{I(x,y,z)=0\}$ is topologically a sphere with
four punctures, and can be seen as the quotient of $\TT$ after $\{\pm\,\Id\,\}$
(i.e.\ one identifies $t$ with $-t$ on $\TT$). This object still admits an
affine extension, though only with a finite translation part (isomorphic
with Klein's 4-group). This affine group in turn is isomorphic with the
group of {\em all} polynomial mappings of 3-space that fix $I(x,y,z)$ --
which gives this simple symmetry analysis a nice application in the context
of nonlinear discrete dynamical systems.

\section*{Summary and Outlook}
The structure of $Gl(2,\ZZ)$ and its relation to algebraic number theory
allowed for a complete classification of (reversing) symmetry groups of
toral automorphisms, while the corresponding investigation of its projective
counterpart, $PGl(2,\ZZ)$, gave the answer for the (reversing) symmetries
of trace maps. It is, at this point, an obvious question how far this
(somehow exceptionally lucky) situation can be generalized.

The correspondence between integer matrices and algebraic number theory
is, of course, not bound to $2\!\times\!2$-matrices, so one can hope
to get some results on $Gl(n,\ZZ)$. This is indeed possible for all matrices
with {\em irreducible} characteristic polynomial, employing again Dirichlet's
unit theorem \cite{BR3}. The cases with reducible polynomials, however, become
increasingly nasty, and their symmetry groups are closely related to the
crystallographic point groups in higher dimensions -- so one should not
expect a simple and general classification here.

As to the question of reversing symmetries,  the perspective is rather negative 
in the sense that reversibility will become more and more an exception with
increasing dimension, because it implies that the spectrum is self-reciprocal
(i.e.\ with $\lambda$ also $1/\lambda$ belongs to the spectrum).  The latter
is usually not the case for the set of algebraic conjugates of an algebraic
integer which constitute the possible spectra of unimodular matrices.
{}Further results along these lines will be reported separately.  

\section*{Acknowledgments}
It is my pleasure to thank John Roberts for his cooperation, and the
German Science Foundation (DFG) for financial support in form of a 
Heisenberg Fellowship.

\end{document}